\documentclass[12pt]{amsart}
\usepackage{amsfonts,amsmath,amssymb,amsthm,verbatim}
\usepackage[latin1]{inputenc}
\theoremstyle{plain} 
\newtheorem{thm}{Theorem}

\newtheorem{cor}[thm]{Corollary}
\newtheorem{lem}[thm]{Lemma}

\theoremstyle{definition}

\newtheorem{ex}{Example}

\theoremstyle{remark}
\newtheorem{remark}{Remark}

\newcommand{\thmref}[1]{Theorem~\ref{#1}}

\newcommand{\lemref}[1]{Lemma~\ref{#1}}

\newcommand{\x}{\mathbf{x}}

\newcommand{\e}{\epsilon}

\ifx\undefined\xspace \usepackage{xspace} \fi
\newcommand{\qr}[1]{\eqref{#1}}

\newcommand{\st}{\mathrel{\,:\,}}
\newcommand{\ie}{{\it i.e.}\xspace}

 \newcommand{\RR}{{\mathbb R}}
\newcommand{\ZZ}{{\mathbb Z}}

\newcommand{\ifemptythenelse}[3]{%
\def\epty{} %
\def\atempempty{#1} %
\ifx\atempempty\epty %
#2 %
\else %
#3 %
\fi
}
\newcommand{\ifnotempty}[2]{\ifemptythenelse{#1}{}{#2}}
\newcommand{\eptypar}[1]{\ifnotempty{#1}{\parens{#1}}}
\newcommand{\gparens}[3]{{\left#1 #2 \right#3}}
\newcommand{\parens}[1]{\gparens({#1})}
\newcommand{\brackets}[1]{\gparens\{{#1}\}}

\newcommand*{\Setof}[1]{\,\brackets{#1}\,}

\renewcommand{\d}{\,{d}}

\newcommand{\Ordo}[1]{{O\eptypar{#1}}}

\numberwithin{equation}{section}
\numberwithin{thm}{section}

\newsavebox{\Chai}
\sbox{\Chai}{\raisebox{0.3ex}[1.7ex][0.75ex]{{$\chi$}}}

\newcommand{\CB}{\mathcal{B}}
\newcommand{\CM}{\mathcal{M}}

\renewcommand{\x}{x}
\newcommand{\CF}{\mathcal{F}}
\newcommand{\CP}{\mathcal{P}}

\newcommand{\var}{\operatorname{var}}

\let\tl=\tilde
\let\mc=\mathcal
\let\X=X
\let\x=x
\let\T=T
\let\e=\hat

\newlength{\ppph}
\newcommand{\px}{\rule{0pt}{\ppph}}

\newcommand{\ii}[1]{^{(#1)}}

\renewcommand{\lll}{\ll^{\tiny\mathrm{loc}}}
\renewcommand{\ae}{almost everywhere\xspace}
\newcommand{\as}{almost surely\xspace}

\def\ppp#1\par{\par\paragraph{#1}}
\begin{document}
\title[countable state shifts and $g$-measures]{Countable state shifts
and uniqueness of $g$-measures
}
\author{Anders Johansson, Anders \"Oberg and Mark Pollicott}
\thanks{A.J. was supported in part by Ume{\aa} University.} 
\thanks{M.P. was supported by an E.U.-Marie Curie chair and a 
Leverhulme Overseas Fellowship.}

\address{Anders Johansson\\ Division of Mathematics and Statistics\\
  University of G\"avle\\ SE-801 76 G\"avle\\ Sweden} \email{ajj@hig.se}

\address{Anders \"Oberg\\
  Department of Mathematics\\ Uppsala University\\ P.O. Box 480 \\ SE-751 06
  Uppsala \\ Sweden} \email{anders@math.uu.se}

\address{Mark Pollicott\\  Mathematics Institute\\ University of Warwick\\
Coventry CV4 7AL\\ UK} \email{mpollic@maths.warwick.ac.uk}


\maketitle
\begin{abstract}
In this paper we present a new approach to studying $g$-measures 
which is based upon local absolute continuity. We extend the result 
in \cite{johob} that square summability of variations of $g$-functions 
ensures uniqueness of $g$-measures. The  first extension is to the 
case of countably many symbols. The second extension is  to some 
cases where $g \geq 0$, relaxing the earlier requirement in \cite{johob} 
that $\inf g>0$.
\end{abstract}
\section{Introduction}\noindent
Let $S$ be a countable discrete set, $X =S^{\ZZ_+}$ the infinite
product with the usual Tychonoff product topology, $\CB$ the
Borel $\sigma$-algebra on $X$ and $T$ the shift on $X$; \ie, $(Tx)_i
= x_{i+1}$ for $i \geq 0$.  The notion of $g$-measures were 
introduced into ergodic theory by Keane in \cite{mK72}.  
We recall that a $g$-function on $X$ is a measurable
function $g$ such that for all $x\in X$,  $\sum_{y\in \T^{-1}x}
g(y) = 1$. For a given $g$-function, a $g$-measure is defined to
be a $\T$-invariant measure $\mu$ in the space of Borel probability
measures $\mc P(X)$, such that $g(x) = {d\mu}/{d\mu \circ T^{-1}}$. 
In the particular case of
$S$ finite and a continuous $g$-function there always exists at
least one $g$-measure. The problem of finding sufficient
conditions on $0 < g < 1$ for which there is a unique $g$-measure
has been extensively studied. In particular, Walters
\cite{walters1} showed that if $\sum_n \var_{n}(\log g) < +\infty$
then there exists a unique $g$-measure. Recently this result was
improved by Johansson and \"Oberg \cite{johob} who showed the same
conclusion   under the weaker assumption $\sum_n (\var_{n}(\log
g))^2 < +\infty$. Their approach was based upon martingale ideas,
rather than the usual transfer operator techniques. Subsequently,
Berger, Hoffman and Sidoravicius \cite{berger} showed the
sharpness of  this result  by adapting a construction of Bramson
and Kalikow \cite{bram}, to show that for any $\epsilon
> 0$ there exist $g$-functions  satisfying $\sum_n (\var_{n}(\log
g))^{2+\epsilon} < +\infty$ and for which there are two distinct
$g$-measures.

In this paper we will consider more general settings. In particular, 
we will allow $S$ to be an infinite set, and we will
also allow the possibility that $g$ takes the value $0$. In this
context, Walters \cite{walters2}, Sarig \cite{sarig1},
\cite{sarig2}, \cite{sarig3}, and Mauldin and Urbanski
\cite{mauldin-urbanski} have proved various existence and
uniqueness results with hyptheses similar to that of summable
variation. (However, in the infinite state case, existence is no
longer automatic due to the lack of compactness of $X^+$.)  In the
present  paper, one of our main results is a  uniqueness result
which, in particular, subsumes the uniqueness result in
\cite{johob}. However, we shall present a new simplified
approach based on local absolute continuity. This method has
been developed by Shiryaev and co-authors (see for instance 
\cite{shiryaev}), in Probability Theory, but seems novel in the 
context of Ergodic Theory.

\section{Predictable ACS criteria}
Consider a discrete filtration $\CF_n\nearrow \CF$, $n\geq 0$, for a
general measure space $(X,\CF)$. For our purposes we may assume that
each $\CF_n$ is countable. 

A measure $\tl{\mu}\in
\CM(X)$ is said to be \emph{locally absolute continuous} with respect
to a second measure $\mu\in \CM(X)$, written $\tl{\mu} \lll \mu$, if
$\tl{\mu}\vert_{\CF_n}\ll \mu\vert_{\CF_n}$, for all $n\geq 0$, where
$\mu\vert_{\CF_n}$ denote the restriction of $\mu$ to the
sub-$\sigma$-algebra $\CF_n$. We write conditional expectation
relative to a probability measure $\mu$ using the integral notation, \ie,
$\int f \d\mu(x\vert \mathcal G)$ denotes the conditional
$\mu$-expectation $d(f\mu\vert_{\mathcal G})/d\mu\vert_{\mathcal G}$ 
of a function $f$ relative to a sub-$\sigma$-algebra
$\mathcal G\subset \CF$.

We will use an approach developed by Shiryaev and his
co-authors for extending local absolute continuity to absolute
continuity for probability measures. To begin, we recall that for
any two probability measures $\mu$ and $\tl\mu$ defined on a
$\sigma$-algebra $\CF$ of a space $X$, the classical Lebesgue
decomposition tells us that we can write $\mu = \lambda_1 +
\lambda_2$, where $\lambda_1, \lambda_2$ are probability measures
on $\CF$ for which $\lambda_1 \ll \tilde \mu$ and $\lambda_2 \perp
\tl \mu$. 

Assuming local absolute continuity allows us to draw a stronger
conclusion.  More precisely, when $\tl\mu \lll \mu$ we can define the
local likelihood ratio process as the $\CF_n$-adapted process given by
$$ Z_n(x) := \frac{d\tl\mu\vert_{\CF_n}}{d\mu\vert_{\CF_n}}. $$ 
Since $Z_n$ is a $\mu$-martingale and a $\tl\mu$-submartingale, the
limit $\lim_n Z_n$ exists $\mu$-\as and $\tl\mu$-\as unless
$Z_\infty := \limsup Z_n$ equals $\infty$. We can therefore write:
$$\tl\mu(B)=\int_B Z_{\infty}\d\mu + \tl\mu(B \cap \{Z_\infty = \infty\}),$$
where $B\in \CF$ (Theorem 4 of \cite{engelbert}; or p.\ 493 of
\cite{shiryaev-book}). From  this decomposition
 we immediately have that:
\begin{equation}\label{zeta}
\tl\mu \ll \mu \iff \tl\mu(Z_{\infty}<\infty)=1
\end{equation}
(cf.\  Theorem 5 of \cite{engelbert}; or p.\ 495 of
\cite{shiryaev-book}). 

For $x\in X$ and $n\geq 1$, let $\alpha_n(x) = Z_n(x)/Z_{n-1}(x)$ and
let
$$ 
d_n(x) := \int (1-\sqrt{\alpha_k (x)})^2\ \d\mu(x\vert\CF_{k-1}). 
$$
The predictable increasing process $B_n(x) = \sum_{k=1}^{n} d_k(x)$,
$n\geq 0$, is referred to as the ``Hellinger process''. 
Kabanov, Lipster and Shiryaev  in
\cite{kabanov} (see also Jacod and Shiryaev \cite{jacod}, p.\ 253,
Theorem 2.36 with $T=\infty$) proved  the following consequence
of \qr{zeta} that they termed a ``predictable ACS-criteria''. 
\begin{lem}\label{lem0}
 If $\tl\mu \lll \mu$, then $\tl\mu \ll \mu$ if and only if $\lim
 B_n(x) <\infty$ with $\tl\mu$-probability one. 
\end{lem}

\newcommand{\shas}{\textit{a.\ s.\xspace}}
\newcommand{\tlog}{\operatorname{tlog}}
\begin{proof}
  We recall the main steps in the simple proof, adapting pages
  496--498 of \cite{shiryaev-book}. By the submartingale property, it
  suffices to show that the log-likelihood process $\log Z_n =
  \sum_{k=1}^n \log \alpha_k$ converges $\tl\mu$-a.s. if only if the
  Hellinger process converges $\tl\mu$-a.s. Furthermore, writing
  $\tlog x = \log x$ if $|\log x| < 1$ and $\operatorname{sign}(\log
  x)$ otherwise, the convergence of $\log Z_n$ occurs precisely when
  the process
$$ Y_n := \sum_{k=1}^n \tlog \alpha_k $$
converges. 

We now claim that $Y_n$ is a $\tl\mu$-submartingale as well. To see this note
first that if $f$ is $\CF_n$-measurable then 
\begin{equation}\label{mula}
  \int f(x) \d\tl\mu(x\vert \CF_{n-1}) = 
  \int \alpha_n(x) f(x) \d\mu(x\vert\CF_{n-1}). 
\end{equation}
Hence, Jensen's inequality gives 
\begin{equation*}
  \int (Y_n - Y_{n-1}) \d\tl\mu(x\vert \CF_{n-1}) 
  = \int \alpha_n \tlog\alpha_n \d\mu(x\vert \CF_{n-1}) \geq 1 \tlog 1
  = 0,
\end{equation*}
since $x\tlog x$ is a convex function for $x\geq 0$. 

Thus $Y_n$ is a submartingale with bounded increments $|Y_n -
Y_{n-1}|<1$. A Doob decomposition $Y_n = A_n + M_n$, where $A_n$ is
predictable and increasing and $M_n$ is a martingale, then
shows that $\lim Y_n < \infty$ $\tl\mu$-a.s. if and only if 
\begin{multline}\label{ffs}
  \sum_n \int \tlog \alpha_n \d\tl\mu(x\vert \CF_{n-1}) + 
  \sum_n \int (\tlog \alpha_n)^2 \d\tl\mu(x\vert \CF_{n-1}) \\=
  \sum_n \int 
  \left[\alpha_n \tlog \alpha_n + \alpha_n (\tlog \alpha_n)^2 \right]
  \d\mu(x\vert \CF_{n-1}) < \infty,
\end{multline}
$\tl\mu$-a.s. (The first sum is $A_n$ and the second sum 
bounds the variance of $M_n$.) 

Finally, we have for $x\geq 0$ that
$$
 (1/C) (1-\sqrt{x})^2 \leq x\tlog x + x(\tlog x)^2 +1-x \leq C (1-\sqrt{x})^2,
$$
for some $C>0$. Since $\int (1-\alpha_n(x)) \d\mu(x|\CF_{n-1}) =
0$, we see that the sum in \qr{ffs} must converge exactly when the
Hellinger process $B_n(x)$ does. 
\end{proof}

For $x\in X$, $[x]_{\CF_n}$, denotes
the minimal element of $\CF_n$ containing $x$. 
We assume that the mapping $[x]_{\CF_k} \to [x]_{\CF_{k-1}}$ has
countable fibers $S_k(x) = \Setof{[y]_{\CF_k}\st [y]_{\CF_k} \subset
  [x]_{\CF_{k-1}}}$ and we denote by  
$p_k$ and $\tl p_k$ the conditional probabilites on $S_k(x)$
induced by $\mu$ and $\tl\mu$, respectively; \ie
\begin{equation}\label{pkdef}
  p_k(x) := \frac{d\mu\vert_{\CF_k}}{\d\mu\vert_{\CF_{k-1}}}(x) = \mu\{
[x]_{\CF_k} \mid [x]_{\CF_{k-1}}\},
\end{equation}
and $\tl p_k$ is defined analogously. 

We may then write
\begin{equation}\label{helld}
    d_n(x) = \rho^2_H(p_k,\tl p_k) := 
    \sum_{[y]\in S_n(x)} (\sqrt{p_k(y)} - \sqrt{\tl p_k(y)})^2, 
\end{equation}
which says that $d_n(x)$ is the squared Kakutani--Hellinger distance between the
probabilities $p_k$ and $\tl p_k$ on the countable set $S_n(x)$. 
To see \qr{helld}, note that $\alpha_k(x) = \tl p_k(x) / p_k(x)$
and, since $\alpha_k$ is $\CF_k$-measurable, we obtain 
\begin{align*}
  \int (1-\sqrt{\alpha_k})^2\ \d\mu(x\vert\CF_{k-1}) &=
  \sum_{[y]\in S_n(x)} 
   p_k(y)(1-\sqrt{\tl p_k(y) / p_k(y)})^2  \\
 &=  \sum_{[y]\in S_n(x)} 
(\sqrt{\tl p_k(y)} - \sqrt{p_k(y)})^2
\end{align*}

\section{Variations of the $g$-function and absolute continuity of
  $g$-chains}

\newcommand{\svar}{\operatorname{svar}}
\newcommand{\Ti}[1]{T^{-#1}} 
In this section we use the predictable
ACS-criteria to derive a criterion for absolute continuity of
$g$-measures on the ``forward algebra'' which we subsequently use to
give sufficient conditions for uniqueness. It is convenient to work
with the natural extension to the two-sided shift $T$ on the space
$X=S^\ZZ$ instead of $X^+ = S^{\ZZ_+}$. Thus $(Tx)_i = x_{i+1}$, $i\in
\ZZ$ and $(\Ti1x)_i = x_{i-1}$.  The one-sided shift is recaptured by
taking the projection $x\to x_+ = (x_0,x_1,x_2,\dots)$ so that $T x_+
= [Tx]_+$. 

Let $\CF^+$ be the $\sigma$-algebra generated by $x_+$, \ie, $\CF^+ =
\lim_n \CF_n^+$, where $\CF_n^+ = \Setof{[x_0, x_1,\dots,x_{n-1}]}$ is
the filtration of \emph{backward} finite cylinders.  A $g$-function is
a $\CF^+$-measurable function $g(x)=g(x_0,x_1,\dots)$ satisfying
$\sum_{ y_+ \in \Ti1(x_+) } g(y) = 1$ for all $x$.  
A $g$-measure is a $T$-invariant probability measure on $X$ such that
$g = d\mu\vert_{\CF^+}/d\mu\vert_{\Ti1\CF^+}$. Equivalently,
$$ g(x) = \lim_n \frac{\mu[x_0,\ldots,x_n]}{\mu[x_1,\ldots,x_n]}, $$
$\mu$-\ae. 
It is clear that any $g$-measure on $\CB(X^+)$ extends
uniquely to a $g$-measure on $X$. 

If we consider the process 
$$
x\ii k := [T^{-k} x]_+ = (x_{-k},x_{-k+1},\dots,x_0,x_1,\dots), \quad
k\geq 0, 
$$ 
we can note that, since $x\ii {n-1} = T x\ii n$, this is a Markov
chain on $X^+$, regardless of the underlying probability measure
$\mu\in\CP(X)$. In this picture, we add the symbol $x_{-n}\in S$
occuring at position $-n$ at ``time'' $n$ so that time runs in the
reverse with the index of the symbol-sequences.  The \emph{initial
  condition} is the distribution of $x\ii0$, \ie $\mu\vert_{\CF^+}$. 
For a $g$-function $g$, we say that a probability $\mu\in\CP(X)$ is a
\emph{$g$-chain} if the transition probabilities $\mu(x\ii
k|x\ii{k-1})$ are given by $g(x\ii k)$, \ie if $\d\mu\vert_{T^k
  \CF^+}/\d\mu\vert_{T^{k-1} \CF^+}= g \circ T^{-k}$, for $k\geq 1$. 
A $g$-measure $\mu$ corresponds to a \emph{stationary} $g$-chain on
$X^+$. 

We want to apply \thmref{lem0} to the filtration $$\CF_n^{-} =
\Setof{[x_{-n}, x_{-n+1},\dots, x_{-1}]} \nearrow \CF^-$$ of finite
\emph{forward} cylinders in $\CB(X)$ and thus derive absolute
continuity of two $g$-chains with respect to the $\sigma$-algebra
$\CF^- = \lim \CF_n^-$ generated by the symbols added during the
evolution of the Markov chain $x\ii n$. 
Recall that for $f: X^+ \to \RR$, we denote
$$\var_{n} f(x) = 
\sup \left\{|f(x) - f(y)| \hbox{ : } x_i = y_i \text{ for } i \leq
  n\right\}.$$ and we also define the {\it s-variation} by
$$ \svar_n f  = \sup_x \left( \sum_{\sigma\in S}
    (\var_{n+1} f(\sigma,x))^2\right)^{1/2}. $$
\begin{lem}\label{lem1} Assume that the $g$-function $g$ satisfies
  $$ \sum_{n=0}^\infty (\svar_n \sqrt{g})^2 < \infty. $$
  Then for any two $g$-chains $\mu$ and $\tl\mu$, 
$\tl\mu\vert_{\CF^-} \ll \mu\vert_{\CF^-}$ provided 
  $\tl\mu\vert_{\CF^-_n} \ll \mu\vert_{\CF^-_n}$, for each $n$. 
\end{lem}
Note that the property $\tl\mu\vert_{\CF^-} \ll \mu\vert_{\CF^-}$ can be
given the following interpretation: There is no test that, based on
observations of the symbols $x_{-1},x_{-2},\dots$ added over time, can
discern with probability one between the two given initial
conditions $\tl\mu\vert_{\CF^+}$ and $\mu\vert_{\CF^+}$. 
\begin{proof}
  Given two $g$-chains $\mu$ and $\tl\mu$ we consider
  the filtration $\CF_n = \CF^-_n$. Our aim is to show that $d_n(x) \leq
  (\svar_n\sqrt g)^2$ since the lemma then follows from \lemref{lem0}. 
  
  We see that the probability
  $p_n(y)$ (or $\tl p_n(y)$), for $[y]_{\CF_n} \subset [x]_{\CF_{n-1}}$,
  in \qr{pkdef} induced from $\mu$ (or $\tl\mu$) is the probability
  $$
    \pi_n(\sigma|x) =
    \mu([\sigma,x_{-n+1},\dots,x_{-1}])/\mu([x_{-n+1},\dots,x_{-1}])
  $$
  of adding the symbol $\sigma=y_{-n}$ at place $-n$ given the symbols
  $x_{-n+1},\dots,x_{-1}$ at places $-n+1,\dots,-1$. 

  Since $\pi_n(\sigma|x) = \int g(x\ii{n}) \d\mu(x\vert[x]_{\CF^+_{n-1}})$
  it is clear that $\pi_n(\cdot|x)$ and
  $\tl\pi_n(\cdot|x)$ are weighted averages of probabilites of the form
  $g(\cdot,y)$, $y\in X^+$, and where $y$ coincide with
  $x_{-n+1},\dots,x_{-1}$ in the first $n-1$ coordinates. Hence, the
  local Hellinger distance $d_n(x)$ satisfies
  $$ 
  d_n(x) = \rho^2_H(\pi_n,\tl\pi_n)
 \leq \sup_{(y,\tl y)} \rho^2_H( g(\cdot,y), g(\cdot,\tl y) ), 
  $$
  with the supremum taken over all pairs $(y,\tl y)$ in $X^+\times
  X^+$ that coincide in the first $n-1$ coordinates. But this supremum
  is clearly less than $(\svar_n \sqrt g)^2$. 
\end{proof}

\section{Conditions for uniqueness of $g$-measures}

We are now ready to consider the problem of uniqueness of
$g$-measures.  Our main result is the following theorem.  
\begin{thm}\label{main}
 Suppose that $g\geq 0$ and that
$\sum_{n=1}^{\infty}\left(\svar_n\sqrt g\right)^2<\infty$. Then any
two distinct ergodic $g$-measures must be  locally incomparable. 
\end{thm}
Local incomparability of $\mu,\tl\mu\in\CM(X)$, means that neither $\mu \lll
\tl\mu$ nor $\tl\mu \lll \mu$.  Local incomparability holds for example when
the transition operator given by $g$ is periodic. 
\begin{proof}
First of all we note that if $\tl\mu\vert_{\CF^-} \ll
\mu\vert_{\CF^-}$ then it follows from the Birkhoff Ergodic Theorem 
that in fact $\tl\mu \ll \mu$; if we are given a density $f =
d\tl\mu\vert_{\CF^-}/d\mu\vert_{\CF^-} \in L^1(\mu)$ then $h =
\lim_{N\to\infty} \frac 1N \sum_{n=0}^{N-1} f\circ T^n$ must be the
density $d\tl\mu/d\mu$, since it follows that $\int_C h \d\mu =
\tl\mu(C)$ for any finite cylinder $C$. 

Secondly, it is well known that no two ergodic measures on a compact
set can be comparable, so we can not have $\mu \ll \tl\mu$ or
$\tl\mu\ll\mu$. This result can used in our context as follows: If $S$
is not finite, then we can denote by $\bar S = S \cup \{\infty\}$ the
one point compactification and then $X$ is contained in the compact
space $\bar X = {\bar S}^{\ZZ}$. The shift $\bar T: \bar X \to \bar X$
is continuous and $\bar X - X$ is an invariant set consisting of
sequences containing the symbol $\infty$.  The ergodic measures on $X$
correspond to ergodic measures on $\bar X$ with $\mu(X) > 0$ (and by
ergodicity $\mu(X) = 1$).  
\end{proof}

If the zero set $\{x \in X \hbox{ : } g(x) = 0\}$  has empty
interior, then every $g$-measure $\mu$ must assign a positive
probability to each cylinder $[x_0,\dots,x_n]$. Hence all
$g$-measures are locally comparable and the above theorem applies. 
In particular, if $g$ is strictly positive then we can deduce
the following corollary. 
\begin{cor}\label{cor1}
  If, in addition to the assumptions of the theorem, $g>0$ then there
  is at most one $g$-measure. 
\end{cor}

\begin{remark}
If  we consider a more general subshift of finite type, then we
would need to impose suitable recurrence conditions on the
associated transition matrix \cite{mauldin-urbanski},
\cite{sarig2}. 
\end{remark}

In the case $S$ is finite and $g$ continuous the last corollary
reduces to the condition of square summability of variations in
\cite{johob}:
\begin{cor}
If $S$ is finite and  $g>0$ satisfies $\sum_n (\var_{n}(\log g))^2
< +\infty$   then there is precisely one $g$-measure. 
\end{cor}
To see this it is enough to note that the compactness of $X$ implies
that $g(x)$ is bounded away from
$0$ and $1$ and that in this case $\var_n \log g$ and $(\svar_{n-1}\sqrt g)^2$ are of the same order.

\begin{ex}\label{ex11}
Let $X^+ = \ZZ_+^{\ZZ_+}$. Fix a sequence $p_i>0$,
$i\geq 1$, such that $\sum_{i=1}^\infty p_i = 1$ and a number $\alpha>0$. Define
$$ g(i,x) = \begin{cases} p_i b(x) & i \geq 1 \\
                          1 - b(x) & i = 0 \end{cases}
$$
where  
$$
   b(x) = \frac1{\zeta(3+\alpha)} \, \sum_{k=1}^\infty  \frac
   1{k^{3+\alpha}}\frac 1{1+x_{k-1}}. 
$$
Then $(\svar_n \sqrt g)^2$
$$
 = \sup_{\substack{\tl x_k = x_k \\ 0\leq k\leq n-1}}
\left\{ \left(\sqrt{b(\tl x)} -
\sqrt{b(x)}\right)^2  +  \left(\sqrt{1-b(\tl x)} -
\sqrt{1-b(x)}\right)^2 \right\} $$
$$
\leq \sqrt{\var_n b(x)} = \Ordo{ n^{1+\alpha/2} }  
$$
and, since this is summable, there is at most one $g$-measure by
Corollary \ref{cor1}. 

However, since 
$$ \lim_{N\to\infty} \left(\limsup_{\substack{x_k\to\infty\\ 0\leq k \leq N}}
  b(x)\right) = 0$$ 
we see that 
$\var_n \log g = \var_n \log b = \infty$, for all $n\geq 1$, 
so any condition on the variations of $\log g$ does not apply
in this case. 
\end{ex}

As we will see in the next two examples, conditions on the variations 
of $g$ are not always necessary to ensure uniqueness of a 
$g$-measure. 
\begin{ex}
  Let $a_n \to a$, with $0 < a_n, a < 1$.  Let $X^+ = \{0,1\}^{\mathbb
    N}$ then we can define $g: X \to (0,1)$ by
\begin{align*}
g(x) &=  a_n &&\text{ if } x_0 = \cdots = x_{n-1} = 0 \text{ and } x_{n} =1\\
&=1- a_n &&\text{ if } x_1 = \cdots = x_{n-1} = 0 \text{ and } x_0
= x_{n} = 1
\end{align*}
Hulse \cite{hulse} observed that there is a unique $g$-measure
without any hypotheses on $a_n$.  In particular, no conditions on
the variations are required. 
\end{ex}

It should perhaps also be noted that under the condition $\sum_{n}
(\svar_n \sqrt g)^2$ $< \infty$ we obtain absolute continuity on the
forward algebra $\CF^{-}$ also for any non-stationary $g$-chain.  It
could aslo be noted that there are $g$-functions that have a unique
$g$-measure, for which some initial conditions gives chains which are
not a.\ c.\ on $\CF^-$.
\begin{ex} Let $S=\{+1,-1\}$, $X=S^\ZZ$ and 
$$ g(\pm1,x_1,x_2,\dots) := \phi(\pm \sum_{i=1}^\infty a_i x_i ) $$
where $\phi(x) = e^x/(e^x + e^{-x})$ and where we require that $a_n>0$,
$\sum_n a_n < \infty$ but
\begin{equation}
  \label{varr}
  \sum_n (\sum_{i=n}^\infty a_i)^2 = \infty. 
\end{equation}
Then it is not too hard to see that if $\mu$ is the distribution of
the Markov chain $x\ii k$, $k\geq 0$, given that $x\ii0 =
(+1,+1,+1,\dots)$ and $\tl\mu$ the distribution starting in $x\ii0 =
(-1,-1,-1,\dots)$, then, by \qr{varr}, $d_n(x)\geq c \cdot
(\sum_{i=n}^\infty a_i)^2$. Hence, by \lemref{lem0} $\mu$ and $\tl\mu$
must indeed be mutually singular on $\CF^-$. However the condition of
Dobrushin (see, for instance, \cite{fernandez}) can be used to deduce
a unique ergodic measure if we in addition assume that $\sum_i a_i <
1$. 
\end{ex}

\section{Existence of $g$-measures}\noindent
In the case of countable shift state spaces, we cannot rely
on the Schauder--Tychonoff fixed point theorem to produce a
$g$-measure, due to the lack of compactness of $X$. 
\begin{ex}  Consider a subshift of finite type $\Sigma \subset X^+$ of
sequences $(x_n)_{n=0}^\infty \in {\mathbb N}^{\mathbb Z^+}$ where 
$x_n$ can be followed by $x_{n+1}$ iff  $|x_n-x_{n+1}|\leq 1$.   The
function $g(x)=1/3$ is a $g$-function. However, the  properties of
the simple random walk on $\mathbb Z$ implies that there is no finite
$g$-measure. 
\end{ex}

In this paper we do not want to rely upon summability of variations
(or local H\"older continuity), as in Sarig \cite{sarig1}, Sarig
\cite{sarig3} and Mauldin and Urbanski \cite{mauldin-urbanski}, to
derive the existence of $g$-measures.  We can do without this
assumption if we assume instead that the $g$-function $g$ can be
continuously extended to a compactification of $X$.  

We give below a
sufficient condition for existence, which is weak enough to
demonstrate that our uniqueness conditions are not vacuous in the case
of countable state shifts.  Let $\bar S$ be a one-point
compactification of $S$, i.e., $\bar S = S \cup \{\infty\}$ and let
$\bar X$ = ${\bar S}^{\ZZ_+}$. It is clear that a continuous function
$f: X\to \RR$ can be continously extended to $\bar X$ if and only if
the following hold: For every $\epsilon >0$ and $n\in\ZZ_+$ there is
some finite set $B\subset S$ such that $|f(x)-f(y)|<\epsilon$
whenever $x_n,y_n\not\in B$ and $y_i=x_i$ for $i\not=n$.
\begin{thm}
  Suppose that $g$ is continuous $g$-function on $X$ such that $g$ can
  be extended to a continous function on $\bar X$. 
  Suppose further that for every $x\in X$ we have
\begin{equation}\label{tight}
g(\sigma x)\leq K \pi(\sigma),
\end{equation}
where $K\geq 1$ and where $\pi$ is a fixed probability measure on the symbol
set $S$, i.e., $\sum_{\sigma} \pi(\sigma)=1$.  Then there exists at
least one $g$-measure on $X=S^{\mathbb Z_+}$. 
\end{thm}

\begin{proof}
  By the Dominated Convergence Theorem in the context of functions in
  $\ell_1(S)$ it follows from \qr{tight} that the continuous extension
  of $g$ must be a $g$-function on $\bar X$, \ie the extended $g$ has,
  in addition to continuity, the property that $\sum_{\sigma \in \bar
    S} g(\sigma,x) = 1$ for any $x\in\bar X$.  We can therefore assume
  a $g$-measure $\mu\in\CP(\bar X)$ implied by the Schauder-Tychonoff
  fixed point theorem.

  For every $\epsilon>0$ we have a finite set $B_{\epsilon} \subset S$
  such that $ \pi(S\setminus B_{\epsilon})\leq \frac{\epsilon}{K}$ and
  let $B_{\epsilon}^X = \{ x_0 \in B_\epsilon \}$ be the sequences in
  $X=S^{\mathbb Z_+}$ which have their last symbols in $B_{\epsilon}$.
  Then we have that
$$\mu(X\setminus B_{\epsilon}^X)\leq \sum_{\sigma\in X\setminus
B_{\epsilon}} \int_X g(\sigma x)\; d\mu(x)\leq \epsilon. $$ which
immediately implies $\mu(\{ x_0 = \infty \}) = 0$ and hence, by
translation invariance, $\mu(\{\exists\, n\st x_n = \infty \}) =
0$. In other words, $\mu(X)=1$ and $\mu$ corresponds to a
$g$-measure in $\mc P(X)$. 
\end{proof}

\begin{remark} It is easy to check that the $g$-function in Example
  \ref{ex11} satisfies the conditions above. We can hence deduce that
  this $g$-function admits precisely one $g$-measure.  
\end{remark}
\begin{remark}
In particular, if $\var_1 (\log g) < +\infty$ then we can fix $x_0
\in X$ and write $g(\sigma x) \leq e^{\var_1(\log g)}g(\sigma
x_0)$.  Thus, the hypotheses hold with $K = e^{\var_1 (\log g)}$
and $\pi(\sigma) = g(\sigma x_0)$. 
\end{remark}

\end{document}